\documentstyle[12pt]{article}

\begin{document}
\newcommand{\ol }{\overline}
\newcommand{\ul }{\underline }
\newcommand{\ra }{\rightarrow }
\newcommand{\lra }{\longrightarrow }
\newcommand{\ga }{\gamma }
\newcommand{\st }{\stackrel }
\newcommand{\scr }{\scriptsize }
\title{\Large\textbf{ Polynilpotent Multipliers of some Nilpotent Products of Cyclic Groups II }}
\author{\textbf{A. Hokmabadi, B. Mashayekhy \footnote{Correspondence: mashaf@math.um.ac.ir} \ \
and F.Mohammadzadeh  } \\
Department of Mathematics,\\
Ferdowsi University of Mashhad,\\
P. O. Box 1159-91775, Mashhad, Iran}
\date{ }
\maketitle
\begin{abstract}
This article is devoted to present an explicit formula for the $c$th
nilpotent multiplier of nilpotent products of some cyclic groups
$G={\bf {Z}}\stackrel{n_1}{*}{\bf
{Z}}\stackrel{n_2}{*}...\stackrel{n_{t-1}}{*}{\bf
{Z}}\stackrel{n_{t}}{*} {\bf {Z}}_{m_{t+1}}\stackrel{n_{t+1}}{*}{\bf
{Z}}_{m_{t+2}}\stackrel{n_{t+2}}{*}...\stackrel{n_{k}}{*}{\bf{Z}}_{m_{k+1}}$,
where $m_{i+1} | m_i$ for all $t+1 \leq i \leq k$ and $c \geq
n_1\geq n_2\geq ...\geq n_t\geq ...\geq n_{k}$ such that $
(p,m_{t+1})=1$ for all prime $p \leq n_1$. Moreover, we compute the
polynilpotent multiplier of the group $G$ with respect to the
polynilpotent variety ${\mathcal N}_{c_1,c_2,...,c_s}$, where $c_1
\geq n_1.$
\end{abstract}
\textbf{Keywords}: Baer invariant; Variety of groups; Nilpotent
multiplier; Cyclic group.\\
\textbf{Mathematics Subject Classification (2000)}: 20E34; 20E10; 20F18. \\
\newpage
\begin{center}
\hspace{-0.65cm}\textbf{1. Introduction and Motivation}\\
\end{center}

Let a group $G$ be presented as a quotient of a free group $F$ by a
normal subgroup $R$. If $\mathcal{V}$ is the variety of
polynilpotent groups of class row $(c_1,...,c_t)$, $ {\mathcal
N}_{c_1,c_2,...,c_t}$, then the Baer invariant of a group $G$ with
respect to this variety which we call it the polynilpotent
multiplier of $G$ is as follows:
$${\mathcal N}_{c_1,c_2,...,c_t} M(G)=\frac{R \cap \gamma_{c_t+1}\circ
...\circ \gamma_{c_1+1}(F)}{[R,\ _{c_1}F,\
_{c_2}\gamma_{c_1+1}(F),...,\ _{c_t}\gamma_{c_{t-1}+1}\circ ...\circ
\gamma_{c_1+1}(F)]},\ \ (1)$$ where $\gamma_{c_t+1}\circ ...\circ
\gamma_{c_1+1}(F)=\gamma_{c_t+1}(\gamma_{c_{t-1}+1}( ...
(\gamma_{c_1+1}(F))...))$ are the term of iterated lower central
series of $F$ (see [7] for further details). In particular, if $t=1$
and $c_1=c$ then the above notation will be the $c$-nilpotent
multiplier of $G$
$$ {\mathcal N}_{c}M(G)=\frac{R \cap \gamma_{c+1}(F)}{[R,\ _cF]}.$$
The case $c=1$ is the much studied Schur multiplier of $G$ and
denoted by $M(G)$.

Historically, there have been several papers from the beginning of
the twentieth century trying to find some structures for the
well-known notion the Schur multiplier and its varietal
generalization the Baer invariant of some famous products of groups,
such as the direct product, the free product and the nilpotent
product. Determining these Baer invariants of a given group is known
to be very useful for classification of groups into isologism
classes. Also structures of Baer invariants are very essential for
studying varietal capability and covering groups.\

 In 1979, Moghaddam [8] found a formula for the
$c$-nilpotent multiplier of a direct product of two groups where
$c+1$ is a prime number or 4. Also, in 1998, G. Ellis [1] presented
the formula for all $c \geq 1.$ In 1997, Mashayekhy and Moghaddam
[9] presented an explicit formula for the $c$-nilpotent multiplier
of a finite abelian group for every $c\geq1$. In 2006, Mashayekhy
and Parvizi [7] gave an explicit formula for the polynilpotent
multiplier of a finitely generated abelian group.\

In 1992, Gupta and Moghaddam [3] calculated the $c$-nilpotent
multiplier of the nilpotent dihedral group of class $n,G_n \cong
{\bf{Z}}_2 \stackrel{n}{*} {\bf{Z}}_2$.(Note that in 2001 Ellis [2]
remarked that there is a slip in the statement and gave the correct
one.) Also, in 2003, Moghaddam, Mashayekhy and Kayvanfar [10]
extended the previous result and calculated the $c$-nilpotent
multiplier of the $n$th nilpotent product of cyclic groups for
$n=$2,3,4 under some conditions. Moreover, Mashayekhy [5] gave an
implicit formula for the $c$-nilpotent multiplier of a nilpotent
product of cyclic groups. Recently, in 2006, Mashayekhy and Parvizi
[12] presented explicit structures for nilpotent and polynilpotent
multipliers of the $n$th nilpotent product of infinite cyclic
groups. Also, the authors [6] obtained an explicit formula for
nilpotent and polynilpotent multipliers of the $n$th nilpotent
product of cyclic groups under some conditions.\

Now in this article, we extend the recent result and first give an
explicit formula for the $c$-nilpotent multiplier of nilpotent
products of some cyclic groups $G=
{\bf{Z}}\stackrel{n_1}{*}...\stackrel{n_{t-1}}{*}{\bf{Z}}\stackrel{n_{t}}{*}
{\bf{Z}}_{m_{t+1}}\stackrel{n_{t+1}}{*}...\stackrel{n_{k}}{*}{\bf{Z}}_{m_{k+1}}$,
where $m_{i+1} \mid m_i $ for $t+1 \leq i \leq k$ and $c \geq n_1
\geq n_2 \geq ...\geq n_{ k} $ such that $(p,m_{t+1})=1$ for all
prime $p \leq n_1$. Second, we compute the polynilpotent
multiplier of the group $G$ with respect to the polynilpotent
variety ${\mathcal N}_{c_1,c_2,...,c_s}$, where  $c\geq n_1$. Note
that by a group $ G = A_1 \stackrel{n_1}{*} A_2 \stackrel{n_2}{*}
... \stackrel{n_l}{*} A_{l+1} $ we mean the group $ G = (...((A_1
\stackrel{n_1}{*} A_2) \stackrel{n_2}{*} A_3) ...
\stackrel{n_l}{*} A_{l+1})$.

\begin{center}
\hspace{-0.65cm}\textbf{ 2. Notation and Preliminaries}\\
\end{center}

Let $ G_{i}= \langle x_{i}| x_{i}^{k_{i}} \rangle \cong {\bf{Z}} _{k
_{i}} $ be the cyclic group of order $ k _{i} $ if $ k _{i} \geq 1
$, and the infinite cyclic group if $ k _{i} = 0 $. The $ n $th
nilpotent product of the family $ \{ G_{i}\} _{i \in I } $ is
defined as follows: $$ \prod ^{\stackrel{n}{*}} _{i \in I }G_{i} =
\frac{ \prod ^{* } _{i \in I } G_{i}}{ \gamma _{n+1}(\prod  ^{* }
_{i \in I }G_{i}) \cap [G_{i}] ^{*} _{i \in I } },$$ where $ \prod
^{* } _{i \in I }G_{i} $ is the free product of the family $ \{
G_{i} \}_{ i \in I } $, and $$ [G_{i}] ^{*} _{i \in I } = \langle [
G_{i}, G_{j}] | i,j \in I, i \neq j \rangle ^{ \prod ^{* } _{i \in I
}G_{i}} $$ is the cartesian subgroup of the free product $ \prod ^{*
} _{i \in I }G_{i} $ which is the kernel of the natural homomorphism
from $ \prod  ^{* } _{i \in I }G_{i} $ to the direct product $ \prod
^{\times } _{i \in I }G_{i} $. For further properties of the above
notation see H. Neuman [12]. Since $ G_{i}^{,}$s are cyclic, it is
easy to see that $ \gamma _{n+1}(\prod ^{* } _{i \in I } G_{i})
\subseteq [G_{i}] ^{*}  $ and hence $ \prod ^{\stackrel{n}{*}} _{i
\in I }G_{i} = \prod ^{* } _{i \in I } G_{i} / \gamma _{n+1}(\prod
^{* } _{i \in I }G_{i}).$

The proof of the main results of this paper rely only on commutator
calculus, so we state some definitions and lemmas
about commutators.\\
\textbf{Definition 2.1} ([4]). The basic commutators on letters
$x_{1}, x_{2}, ...,x_{n},...$ are defined
as follows: \\
(i) The letters $x_{1}, x_{2}, ...,x_{n}, ...,$ are basic
commutators of weight one, ordered by setting $x_{i} < x_{j}$, if
$i<j$.\\
(ii) If basic commutators $c_{i}$ of weight $w(c_{i})< k$ are
defined and ordered, then define basic commutators of weight $ k $
by the following rules. $[c_{i}, c_{j}]$ is a basic commutator of
weight $ k $  if and only if\\
1. $w(c_{i}) + w(c_{j})=k ;$\\
2. $c_{i} > c_{j};$\\
3. if $c_{i}=[c_{s}, c_{t}]$, then $c_{j} \geq c_{t}$.\\
Then we will continue the order by setting $c \geq c_{i}$, whenever
$w(c) \geq w(c_{i}),$ fixing any order among those of
weight $ k $, and finally numbering them in order.\\
\textbf{Theorem 2.2} ([4]). Let $F$ be the free group on $\{x_{1},
x_{2}, ...,x_{d}\}$, then for all $1 \leq i \leq n$, $$
\frac{\gamma_{n}(F)}{\gamma_{n+i}(F)} $$ is a free abelian group
freely generated by the basic commutators of weights $ n, n+1,...,
n+i-1 $  on the letters
$\{x_{1}, x_{2}, ...,x_{d}\}$.\\
\textbf{Theorem 2.3} (Witt formula [4]). The number of basic
commutators of weight $n$  on $ d $  generators is given by the
following formula $$ \chi_{n}(d)= \frac{1}{n}\sum_{m|n}\mu(m)
d^{\frac{n}{m}},$$ where $\mu(m)$ is the M\"{o}bius function, which
defined to be
 $$\mu(m)=\left\{\begin{array}{ll}
 1& ; m=1,\\
 0& ; m=p_1^{\alpha_1}...p_k^{\alpha_k} \ \ \  \exists \alpha_{i} > 1,\\
(-1)^s& ; m=p_1...p_s, \end{array}\right.$$ where the $p_i$ are
distinct prime numbers.\\
\textbf{Lemma 2.4}. Let $1\leq w < r$ and $(p,r)=1$ for all prime
$p$
less than or equal to $w$, then $r$ divides ${r \choose w}$.\\
\textbf{Proof}. Clearly ${r \choose w}=r( \frac{(r-1)...(r-w+1)}{1
\times 2 \times...\times w})$ is an integer. For any prime $p\leq
w$, since $p \mid \!\!\!\!/r$ ,  $p |(r-1)...(r-w+1)$. Thus $1
\times 2 \times ... \times w |(r-1)...(r-w+1)$ and hence the result
holds.$\Box$

The following consequences of the collecting process are very
vital in the proof of the main result.\\
\textbf{Lemma 2.5} (R. R. Struik [13]). Let $x,y$ be any elements of
a group and let $c_1,c_2,...$ , be the sequence of basic commutators
of weight at least two in $x$ and $[x,y]$, in ascending order. Then
$$[x^{n},y]=[x,y]^{n}c_1^{f_1(n)}c_2^{f_2(n)}...c_i^{f_i(n)}...,\
\ \ \ (2)$$where$$f_i(n)=a_1{n \choose 1}+ a_2{n \choose
2}+...+a_{w_i}{n \choose w_i}, \ \ \ \ ( 3 ) \ $$with $a_i \in
\bf{Z}$, and $w_i$ is the weight of $c_i$ in $x$ and $[x,y]$. If the
group is nilpotent, then the expression in (2) gives an identity,
and the sequence of commutators terminates. \\
\textbf{Lemma 2.6} (R. R. Struik [13]). Let $\alpha$ be a fixed
integer, and $G$ be a nilpotent group of class at most $n$. If $b_j
\in G$ and $r<n$, then
$$[b_1,..,b_{i-1},b_i^{\alpha},b_{i+1},...,b_r]=[b_1,...,b_r]^{\alpha}c_1^{f_1(\alpha)}c_2^{f_2(\alpha)}...,$$
where the $c_k$ are commutators in $b_1,...,b_r$ of weight strictly
greater than $r$, and every $b_j$, $1\leq j \leq r$ appears in each
commutator $c_k$, the $c_k$ listed in ascending order. The $f_i$ are
of the form (3), with $a_j \in \bf{Z}$, and
$w_i$ is the weight of $c_i$ (in the $b_i$) minus $(r-1)$.\\

\begin{center}
\hspace{-0.65cm}\textbf{ 3. Main Results}
\end{center}

Let $ A_{i}= \langle x_{i}| x_{i}^{m_{i}} \rangle $ be cyclic groups
such that $A_i \cong {\bf{Z}}$ ($m_i=0$) for $ 1 \leq i \leq t$ and
$A_i \cong {\bf{Z}_{m_i}}$ ($m_i>1$) for $ t+1 \leq i \leq k+1$. Let
$$ 1 \rightarrow R _{i}= \langle x_{i}^{m_{i}} \rangle \rightarrow
F_{i}= \langle x_{i} \rangle \rightarrow A_{i} \rightarrow 1 $$ be a
free presentation for $ A_{i}$.

Put $ G = A_1 \stackrel{n_1}{*} A_2 \stackrel{n_2}{*} ...
\stackrel{n_k}{*} A_{k+1} $, where $n_1 \geq n_2 \geq ... \geq n_k$.
We are going to find a free presentation for $G$. We also use the
notation $\bar{F}_s = F_1*F_2...*F_s $ for $2 \leq s
\leq k+1$ and $F=\bar{F}_{k+1} $.\\
\textbf{Lemma 3.1}. Keep the above notation and put $ S= \langle
x_{i}^{m_{i}} | 1 \leq i \leq k+1 \rangle ^F $ and $R=S
\gamma_{n_1+1}(F) \prod^k_{i=2}( \gamma_{n_i+1}(\bar{F}_{i+1}) \cap
[\bar{F}_{i},F_{i+1}])$. Then $ 1 \rightarrow R  \rightarrow
F \rightarrow G \rightarrow 1 $ is a free presentation of $G$.\\
\textbf{Proof}. We use induction on $k$. When $k=1$, it is easy to
see that $G\cong F/S \gamma_{n_1+1}(F)$. Suppose that $k >1$ and
that the result has been proved for the group $A=A_1
\stackrel{n_1}{*} A_2 \stackrel{n_2}{*} ... \stackrel{n_{k-1}}{*}
A_{k}$. Then $A$ has a free presentation as follows $ 1 \rightarrow
R'  \rightarrow \bar{F}_{k} \rightarrow A \rightarrow 1 $, where
$R'=S' \gamma_{n_1+1}(\bar{F}_{k}) \prod^{k-1}_{i=2}(
\gamma_{n_i+1}(\bar{F}_{i+1}) \cap [\bar{F}_{i},F_{i+1}])$ and $ S'=
\langle x_{i}^{m_{i}} | 1 \leq i \leq k \rangle ^{\bar{F}_{k}}$. On
the other hand, there exist the following sequence
$$F= \bar{F}_{k}*F_{k+1} \stackrel{\varphi}{\rightarrow} A*A_{k+1}
\stackrel{\psi}{\rightarrow} A \stackrel{n_k}{*} A_{k+1} \rightarrow
1. $$ It is easy to see that $$ker(\psi \circ \varphi)=\langle
x_{k+1}^{m_{k+1}} , R' , \gamma_{n_k+1}(F) \cap
[\bar{F}_{k},F_{k+1}] \rangle^F .$$ So
$$ker(\psi\circ\varphi)=\langle  S \gamma_{n_1+1}(\bar{F}_{k}) \prod^{k}_{i=2}(
\gamma_{n_i+1}(\bar{F}_{i+1}) \cap [\bar{F}_{i},F_{i+1}]) \rangle^F
.$$ We also note that

\begin{eqnarray*}
\gamma_{n_1+1}(F) &=& \gamma_{n_1+1}(\bar{F}_{k})
(\gamma_{n_1+1}(F)\cap [\bar{F}_{k},F_{k+1}])\\ &\subseteq&
\gamma_{n_1+1}(\bar{F}_{k})(\gamma_{n_k+1}(F)\cap
[\bar{F}_{k},F_{k+1}]).
\end{eqnarray*}
Hence  $$ker(\psi\circ\varphi)=\langle  S \gamma_{n_1+1}(F)
\prod^{k}_{i=2}( \gamma_{n_i+1}(\bar{F}_{i+1}) \cap
[\bar{F}_{i},F_{i+1}]) \rangle^F =R $$ and this completes the
proof.$\Box$\\

Keeping the previous notation put $$
M_i=\gamma_{n_i+1}(\bar{F}_{i+1})\cap [\bar{F}_{i} ,F_{i+1}]$$
$$N_i=\gamma_{n_i+c+1}(\bar{F}_{i+1})\cap [\bar{F}_{i} ,F_{i+1}].$$
Then by the previous lemma the $c$-nilpotent multiplier of $ G = A_1
\stackrel{n_1}{*} A_2 \stackrel{n_2}{*} ... \stackrel{n_k}{*}
A_{k+1} $ is as follows for all $c \geq n_1 $
$$ {\mathcal N}_{c}M( G)= \frac{\gamma _{c+1}(F)}{[S, \ _c F] \gamma_{n_1+c+1}(F)
\prod^k_{i=2}[M_i, \ _c F]}. $$

We now prove the following lemmas to compute the $c$-nilpotent
multiplier of $G$.\\
\textbf{Lemma 3.2}. $[M_i, \ _c\bar{F}_{i+1}]=N_i$.\\
\textbf{Proof}. Clearly $[M_i, \ _c\bar{F}_{i+1}]\subseteq N_i$. For
the reverse inclusion we have
\begin{eqnarray*}
\gamma_{n_i+c+1}(\bar{F}_{i+1})&=&[\gamma_{n_i+1}(\bar{F}_{i+1}) , \
_c\bar{F}_{i+1}]\\&=&[\gamma_{n_i+1}(\bar{F}_{i}), \
_c\bar{F}_{i+1}][M_i, \ _c\bar{F}_{i+1}]\\&=&
\gamma_{n_i+c+1}(\bar{F}_{i}) \prod_{\exists t_j \ s.t. \
t_j=i+1}[\gamma_{n_i+1}(\bar{F}_{i}),F_{t_1},...,F_{t_c}][M_i, \
_c\bar{F}_{i+1}].
\end{eqnarray*}
Moreover
\begin{eqnarray*}
[\gamma_{n_i+1}(\bar{F}_{i})
,F_{t_1},...,F_{t_s},F_{i+1},...,F_{t_c}]&=&[
F_{i+1},[\gamma_{n_i+1}(\bar{F}_{i}),F_{t_1},...,F_{t_s}],...,F_{t_c}]\\&\subseteq
&[F_{i+1},[\gamma_{n_i+1}(\bar{F}_{i+1}), \
 _s\bar{F}_{i+1}],...,\bar{F}_{i+1}]\\& \subseteq &[F_{i+1}, \
_{n_i+s+1}(\bar{F}_{i+1}), \ _{c-s}\bar{F}_{i+1}]\\&\subseteq&[M_i,
\ _c\bar{F}_{i+1}].
\end{eqnarray*}
Therefore
$\gamma_{n_i+c+1}(\bar{F}_{i+1})=\gamma_{n_i+c+1}(\bar{F}_{i})[M_i,
\ _c\bar{F}_{i+1}].$ On the other hand,
$$\gamma_{n_i+c+1}(\bar{F}_{i+1})=\gamma_{n_i+c+1}(\bar{F}_{i})
(\gamma_{n_i+c+1}(\bar{F}_{i+1})\cap [\bar{F}_{i},F_{i+1}]).$$
Hence the result holds.$\Box$\\
\textbf{Lemma 3.3}. With the notation and assumption above we have
$$\prod_{i=2}^{k}[M_i, \ _cF]=\prod_{i=2}^{k}N_i.$$
\textbf{Proof}. The previous lemma implies that $N_i=[M_i, \
_c(\bar{F}_{i+1})] \subseteq [M_i, \ _cF]$ for all $2 \leq i \leq
k$. For the reverse inclusion we have
$$ [M_i, \ _cF]=\prod_{D_j\in L,\\
1\leq j\leq c}[M_i,D_1,...,D_c] \subseteq \prod_{i=1}^{k}N_i,$$
where $L=\{ F_j | 1 \leq j \leq k+1 \}$.
 Note that
$[M_i,D_1,...,D_c]\subseteq N_l$, where $l=max \{i,m\}$ such that
$m=max\{t \mid F_t \in \{D_1,...,D_c\}\}.\Box $

Now the above lemma implies that
\begin{eqnarray*}
{\mathcal N}_{c}M(G)&\cong& \frac{\gamma_{c+1}(F)}{[S, \
_cF]\prod_{i=2}^{k}[M_i, \
_cF]\gamma_{n_1+c+1}(F)}\\&=&\frac{\gamma_{c+1}(F)}{[S, \
_cF]\prod_{i=2}^{k}N_i \gamma_{n_1+c+1}(F)}\\
&\cong& \frac{\gamma_{c+1}(F)/\gamma_{n_1+c+1}(F)}{([S, \
_cF]\prod_{i=2}^{k}N_i \gamma_{n_1+c+1}(F))/\gamma_{n_1+c+1}(F)}.\ \
\ (4)
\end{eqnarray*}
Also by Theorem 2.2 $\gamma_{c+1}(F)/\gamma_{n_1+c+1}(F)$ is a free
abelian group with the basis of all basic commutators of weight
$c+1,...,c+n_1$ on letters $x_1,...,x_{k+1}$. Now we are going to
find a suitable basis for the group
  $$\frac{[S, \ _cF]\prod_{i=2}^{k}N_i
\gamma_{n_1+c+1}(F)}{\gamma_{n_1+c+1}(F)}.$$
\textbf{Lemma 3.4}.
With the previous notation $ (\prod_{i=2}^{k}N_i
\gamma_{n_1+c+1}(F))/\gamma_{n_1+c+1}(F)$ is a free abelian group
with the basis $ \bigcup _{i=2} ^{k}E_i$, where
$$E_i=\{ all \ basic \ commutators \ of \ weight \
c+n_i+1,...,c+n_1$$ $$ \ \ \ \ \ on \ letters \ x_1,...,x_{i+1} \
such \ that \ x_{i+1} \ appears \}$$
\textbf{Proof}. Since $E_i$'s
are distinct, it is enough to show that $E_i$ is a basis for $( N_i
\gamma_{n_1+c+1}(F))/\gamma_{n_1+c+1}(F)$ for all $2 \leq i \leq k
$. We know that $ \gamma_{n_i+c+1}(\bar{F}_{i+1})=
\gamma_{n_i+c+1}(\bar{F}_{i})N_i$ and $\gamma_{n_i+c+1}(\bar{F}_{i})
\cap N_i=1$. On the other hand, $\gamma_{n_i+c+1}(\bar{F}_{i+1})$ is
generated by all basic commutators of weight $c+n_i+1,...,c+n_1$ on
letters $x_1,...,x_{i+1}$ and $\gamma_{n_i+c+1}(\bar{F}_{i})$ is
generated by all basic commutators of weight $c+n_i+1,...,c+n_1$ on
letters $x_1,...,x_{i}$  , modulo $\gamma_{n_1+c+1}(F)$ . So
the result follows.$\Box$\\

We define $\rho_1(S)=S$ and $\rho _{r+1}(S)=[ \rho _{r}(S),F]$,
inductively.\\
\textbf{Lemma 3.5}. With the notation and assumption above, if $( p,
m_{t+1})=1$, for all prime $p$ less than or equal to $n_1+1-i$, then
$(\rho _{c+i}(S)\gamma _{c+n_1+1}(F))/$\\$ (\rho _{c+i+1}(S)\gamma
_{c+n_1+1}(F))$ is a free abelian group with a basis
$D_{i,1}\cup...\cup D_{i,k+1-t}$, where\\ $$D_{i,j}=\{ b^{m_{t+j}} |
 \ b\ is\ a\ basic\ commutator\ of\ weight\ c+i\ on\ $$ $$ \ \ \ \
 \ x_1,...,x_t,...,x_{t+j}\ such\ that \ x_{t+j}\ appears\ in\ b \ \},$$
 for all $1 \leq i \leq n_1$. \\
\textbf{Proof}. Using collecting process (see [4]), it is easy to
see that
$$(\rho _{c+i}(S)\gamma _{c+n_1+1}(F))/ (\rho _{c+i+1}(S)\gamma
_{c+n_1+1}(F))$$ is generated by all basic commutators of weight
$c+i, ...,c+n_1$ on $x_1,...,x_t,..,x_{k+1}$ such that one of the
$x_{t+1}^{m_{t+1}},...,x_{k+1}^{m_{k+1}}$ appears in them. Also it
is routine to check that all the above commutators of weight grater
than $c+i$ belong to $\rho _{c+i+1}(S)$. Now, we are going to show
that if $b'$ is a basic commutator of weight $c+i$ on $
x_1,...,x_t,...,x_{k+1}$ such that $x_{t+j}^{m_{t+j}}$, appears in
it, then $$b'\equiv b^{m_{t+j}} \ \ \ \ \ \pmod{\rho
_{c+i+1}(S)\gamma _{c+n_1+1}(F)}, \ \ (5) $$ where $b$ is a basic
commutator of weight $c+i$ on $ x_1,...,x_t,...,x_{k+1}$ such that
$x_{t+j}$ appears in it. (Note that $b$ is actually a real basic
commutator according to the definition and $b'$ is the same as $b$
but a letter $x_{t+j}$ with exponent $m_{t+j}$.) In order to prove
the above claim, first we use reverse induction on $d$ $(i+1 \leq d
\leq n_1)$ to show that if $u$ is an outer commutator of weight
$c+d$ on $x_1,...,x_t,...,x_{k+1}$ such that $x_{t+j}$ appears in
$u$, then $u^{m_{t+j}}\in \rho _{c+i+1}(S)\pmod{\gamma
_{c+n_1+1}(F)}.\ \ (6)$
\\let $d=n_1$ and $u=[...,x_{t+j},...]$ then clearly $u^{m_{t+j}}
\equiv [...,x_{t+j}^{m_{t+j}},...]\in \rho _{c+i+1}(S)\pmod {\gamma
_{c+n_1+1}(F)}$.\\
Now, suppose the above property holds for every $d>d'$. We will show
that the property (5) holds for $d'$. Let $u=[u_1,u_2]$ be an outer
commutator of weight $c+d'$ where $x_{t+j}$ appears in $u_1$, say.
Then by Lemma 2.5, we have $$u^{m_{t+j}} \equiv
[u_1^{m_{t+j}},u_2]v_1^{f_1(m_{t+j})}...v_h^{f_h(m_{t+j})} \pmod{
\gamma _{c+n_1+1}(F)},$$where $v_s$ is a basic commutator of weight
$w_s$ in $u_1$ and $[u_1,u_2]$. Thus $v_s$ is an outer commutator of
weight greater than $c+d'$ and less than $c+n_1+1$ on $
x_1,...,x_t,...,x_{k+1}$ such that $x_{t+j}$ appears in it. Since
$m_{t+j}|m_{t+1}$, by hypothesis we have $( p, m_{t+j})=1$, for all
prime $p$ less than or equal to $n_1+1-i$. Also, it is easy to see
that $w_s \leq (c+n_1+1)-(c+d'-1)=n_1+1-d'+1 \leq n_1+1-i$.
Therefore by Lemma 2.4 $m_{t+j}|f_s(m_{t+j})$ and so by induction
hypothesis $v_s^{f_s(m_{t+j})} \in \rho _{c+i+1}(S)\pmod{ \gamma
_{c+n_1+1}(F)$}. Hence, by repeating the above process if
$u=[...,x_{t+j},...]$, then $u^{m_{t+j}} \equiv
[...,x_{t+j}^{m_{t+j}},...]v_1'^{f_1'(m_{t+j})}...v
_h'^{f_h'(m_{t+j})}\in \rho _{c+i+1}(S) \pmod{ \gamma _{c+n_1+1}(F)}.$ \\
 Now using (6), Lemma 2.6 and some routine commutator calculus,
the congruence (5) holds. By  Theorem 2.2 distinct basic commutators
are linear independent and so are distinct powers of them. Since
$m_{t+j+1} | m_{t+j}$, for all $1 \leq j \leq k+1-t$, the set $
\bigcup ^{k+1-t}_{j=1}D_{i,j}$ is a generating set for $(\rho
_{c+i}(S)\gamma _{c+n_1+1}(F))/$\\$(  \rho _{c+i+1}(S)\gamma
_{c+n_1+1}(F))$ and
hence is a basis for it.$\Box$\\
\textbf{Lemma 3.6}. If $(p,m_{t+1})=1$ for all prime $p$ less than
or equal to $n_1$, then $(\rho _{c+1}(S)\gamma _{c+n_1+1}(F))/
\gamma _{c+n_1+1}(F)$ is a free abelian group with a basis $ \bigcup
^{k}_{i=t}C_{i}$, where $$C_i=\{ b^{m_{i+1}} |
 \ b\ is\ a\ basic\ commutator\ of\ weight\ c+1,...,c+n_1 $$ $$ \ \ \ on \ letters \
x_1,...,x_{i+1} \ such \ that \  x_{i+1} \  appears \ in \ b \}.$$
\textbf{Proof}. Put $$A_i= \frac{\rho _{c+i}(S)\gamma
_{c+n_1+1}(F)}{ \rho _{c+i+1}(S)\gamma _{c+n_1+1}(F)}, B_i=\frac{
\rho _{c+1}(S)\gamma _{c+n_1+1}(F)}{ \rho _{c+i+1}(S)\gamma
_{c+n_1+1}(F)}.$$ Then clearly the following exact sequence exists
for all $1 \leq i \leq n_1$
$$0 \rightarrow A_i \rightarrow B_i \rightarrow B_{i-1}
\rightarrow 0. $$ By Lemma 3.5 $B_1$ is a free abelian, so the exact
sequence$$0 \rightarrow A_2 \rightarrow B_2 \rightarrow B_1
\rightarrow 0 $$is a split exact sequence, and hence $B_2 \cong A_2
\oplus B_1 $. Also, by Lemma 3.5 every $A_i$ is free abelian, so by
induction every $B_i$ is free abelian and $$\frac{\rho
_{c+1}(S)\gamma _{c+n_1+1}(F)}{\gamma _{c+n_1+1}(F)}=B_{n_1} \cong
A_{n_1} \oplus A_{n_1-1} \oplus...\oplus A_2 \oplus A_1. $$ Now,
according to the basis for $A_i$ presented in Lemma 3.5, the group
$(\rho _{c+1}(S)\gamma _{c+n_1+1}(F))/ \gamma _{c+n_1+1}(F)$ is a
free abelian group with a basis $$ \bigcup _{i=1}^{n_1}( \bigcup
^{k+1-t}_{j=1}D_{i,j})= \bigcup ^{k}_{i=t}C_{i}.  \Box$$

Lemma 3.4 and Lemma 3.6 implies that
$\bigcup _{i=t} ^{k}C_i \cup \bigcup^{k}_{i=2}E_i$ generate the dinominator of 4,
but some elements of $C_i$ which are generated by
elements of $\bigcup^{k}_{i=2}E_i$ can be omitted and other elements are linear independent.
 Therefore we have the following theorem. \\
\textbf{Theorem 3.7}. Keeping the notation of Lemma 3.4 if
$(p,m_{t+1})=1$ for all prime $p$ less than or equal to $n_1$, then
$([S, \ _cF]\prod_{i=2}^{k}N_i
\gamma_{n_1+c+1}(F))/\gamma_{n_1+c+1}(F)$ is a free abelian group
with the basis $(\bigcup _{i=2} ^{k}E_i)\cup (\bigcup_{i=t}
^{k}L_i)$, where  $$L_i=\{ b^{m_{i+1}}| b \ is \ a \ basic \
commutator \ of \ weight \ c+1,...,c+n_i $$ $$ \ \ \ on \ letters \
x_1,...,x_{i+1} \ such \ that \  x_{i+1} \  appears \ in \ b \}.$$

Now we are ready to state and prove the first main result of the
paper.\\
\textbf{Proposition 3.8}. Let $
G=A_1\stackrel{n_1}{*}...\stackrel{n_k}{*}A_{k+1}$ such that
$A_i\cong {\bf Z}$ for $1\leq i \leq t$ and $A_j\cong {\bf Z}_{m_j}$
for  $ t+1 \leq j \leq k+1$. Let $c\geq n_1 \geq ...
\geq n_k$ and $m_{k+1}| m_k|...|m_{t+1}$ and $(p,m_{t+1})=1$ for all prime $p \leq n_1$. Then \\
$${\mathcal N}_{c}M(G)\cong {\bf Z}^{(u)}\oplus {\bf
Z}_{m_{t+1}}^{(f_t)}\oplus...\oplus {\bf Z}_{m_{k+1}}^{(f_k)},$$
where $u = \sum _{j=1} ^{n_{t-1}}\chi _{c+j}(t)+\sum _{i=1} ^{{t-2}}
\sum _{j=n_{i+1}+1} ^{n_{i}}\chi _{c+j}(i+1)$ and
$f_s=\sum _{j=1} ^{n_{s}}(\chi _{c+j}(s+1)-\chi _{c+j}(s))$.\\
\textbf{Proof}. Applying the structure of ${\mathcal N}_{c}M(G)$
given by (4) and Theorem 3.7 and note that the $E_i$ and the $L_i$
are mutually disjoint and
$$|E_i|=\sum _{j=n_i+1} ^{n_1}(\chi _{c+j}(i+1)-\chi _{c+j}(i)),$$
$$|L_i|=\sum _{j=1} ^{n_i}(\chi _{c+j}(i+1)-\chi _{c+j}(i)),$$ we
have

\begin{eqnarray*}
u &=& \sum^{n_1}_{j=1}
\chi_{c+j}(k+1)-(\sum^{k}_{i=t}\sum^{n_i}_{j=1}(
\chi_{c+j}(i+1)-\chi_{c+j}(i))+ \sum^{k}_{i=2}\sum^{n_1}_{j=n_i+1}(
\chi_{c+j}(i+1)-\chi_{c+j}(i)))\\&=& \sum^{n_1}_{j=1}
\chi_{c+j}(k+1)-\sum^{n_1}_{j=1}\sum^{k}_{i=t}(
\chi_{c+j}(i+1)-\chi_{c+j}(i))-\sum^{t-1}_{i=2}\sum^{n_1}_{j=n_i+1}(
\chi_{c+j}(i+1)-\chi_{c+j}(i)))\\&=&\sum^{n_1}_{j=1}
\chi_{c+j}(t)-\sum^{t-1}_{i=2}\sum^{n_1}_{j=n_i+1}(
\chi_{c+j}(i+1)-\chi_{c+j}(i)))\\&=&\sum^{n_1}_{j=1}
\chi_{c+j}(t)-\sum^{t-1}_{i=2}\sum^{n_1}_{j=n_i+1}
\chi_{c+j}(i+1)+\sum^{t-2}_{i=1}\sum^{n_1}_{j=n_i+1}
\chi_{c+j}(i+1)\\&=&\sum _{j=1} ^{n_{t-1}}\chi _{c+j}(t)+\sum _{i=1}
^{{t-2}} \sum _{j=n_{i+1}+1} ^{n_{i}}\chi _{c+j}(i+1)
\end{eqnarray*}
Hence the result follows. $\Box $

In order to state and prove the second main result we apply the
following lemma.\\
\textbf{Lemma 3.9}. Let $G$ be a nilpotent group of class $n \leq
c_1$, and $c_2,...,c_s \geq 1$. Then $${\mathcal N}_{c_1,c_2
...,c_s}M(G)={\mathcal N}_{c_s}M(...{\mathcal N}_{c_2}M({\mathcal
N}_{c_1}M(G))...).$$
\textbf{Proof}. Let $F/R$ be a free
presentation of $G$. Since $\gamma_{c_1+1}(F)\leq \gamma_{n+1}(F)
\leq R$, ${\mathcal N}_{c_1}M(G)=\gamma _{c_1+1}(F)/[R,\ _{c_1}F]$.
Now, we can consider $\gamma _{c_1+1}(F)/[R,\ _{c_1}F]$ as a free
presentation for ${\mathcal N}_{c_1}M(G)$ and hence$${\mathcal
N}_{c_2}M({\mathcal N}_{c_1}M(G))= \frac{\gamma _{c_2+1}(\gamma
_{c_1+1}(F))}{[R,\ _{c_1}F,\ _{c_2}\gamma _{c_1+1}F]}.$$ Therefore
by (1) we have
$${\mathcal N}_{c_1,c_2}M(G)={\mathcal N}_{c_2}M({\mathcal
N}_{c_1}M(G)).$$ By continuing the above process we can conclude the
result.$\Box$

Note that, by a similar method of [9], we can obtain the structure
of the $c$-nilpotent multiplier of a finitely generated abelian
group as the following statement. Also, this is a corollary of
Proposition 3.8 when $n_1=n_2=...=n_k=1$\\
\hspace{-0.65cm}\textbf{Corollary 3.10}. Let $G= {\bf{Z}}^{(m)}
\oplus {\bf{Z}}_{m_{1}} \oplus ...\oplus{\bf{ Z}}_{m_{k}},$ be a
finitely generated abelian group, where $m_{i+1}| m_i$ for all
$1\leq i\leq k-1$. Then\\
${\mathcal N}_{c}M(G)= {\bf{Z}}^{(b_m)} \oplus
{\bf{Z}}_{m_{1}}^{(b_{m+1}-b_m)} \oplus ...\oplus{\bf{ Z}}_{m_{k}}^{
(b_{m+k}-d_{m+k-1})},$ where $b_i=\chi_{c+1}(i)$.

 Now the second main result of the
paper follows immediately by
the above lemma, Proposition 3.8 and repeated usage of Corollary 3.10.\\
\textbf{Proposition 3.11}. Let $G = A_1 \stackrel{n_1}{*} A_2
\stackrel{n_2}{*} ... \stackrel{n_k}{*} A_{k+1}$ such that
$A_i\cong {\bf Z}$ for $1\leq i \leq t$ and $A_j\cong {\bf
Z}_{m_j}$ for  $ t+1 \leq j \leq k+1$. Let $c_1\geq n_1 \geq n_2
\geq ...\geq n_k$ and $m_{k+1}| m_k|...|m_{t+1}$ and
$(p,m_{t+1})=1$ for all prime $p \leq n_1$. Then the structure of
the polynilpotent multiplier of $G$ is
$${\mathcal N}_{c_1,c_2,...,c_s}M(G)\cong {\bf{Z}}^{(d_0)} \oplus {\bf{Z}}_{m_{t+1}}^{(d_{t}-d_o)}
\oplus ...\oplus{\bf{ Z}}_{m_{k+1}}^{ (d_{k}-d_{k-1})},$$  where
$d_i=\chi_{c_s+1}(...(\chi_{c_2+1}(u+\sum ^{i}_{j=t}f_j))...),$ for
all $ t \leq i \leq k$, $d_0=\chi_{c_s+1}(...(\chi_{c_2+1}(u))...)$
and $u$ and $f_i$ are the notation used in Proposition 3.8.


\begin{thebibliography}{19}
\bibitem{1} G. Ellis, On Groups with a Finite Nilpotent Upper
Central Quotient, {\it Arch. Math.} {\bf 70}, (1998) 89-96.
\bibitem{2} G. Ellis, On the Relation Between Upper Central Quotients and Lower Central Series of a
Group, {\it Trance. Amer. Math. Soc.} {\bf 353}, (2001) 4219-4234.
\bibitem{3} N. D. Gupta and M. R. R. Moghaddam, Higher Schur
Multiplicators of Nilpotent Diheadral Groups, {\it C. R. Math. Rep.
Acad. Sci. Canada} XIV {\bf 5}, (1992) 225-230.
\bibitem{4} M. Hall, {\it The Theory of Groups}, Macmillian Company, New York, 1959.
\bibitem{5} B. Mashayekhy, Some Notes on the Baer-invarient of a
Nilpotent Product of Groups, {\it Jornal of Algebra} {\bf 235},
(2001) 15-26.
\bibitem{6} B. Mashayekhy, A. Hokmabadi and F. Mohammadzade, Polynilpotent
Multipliers of some Nilpotent Products of Cyclic Groups, Submitted.
\bibitem{7}  B. Mashayekhy and M. Parvizi, Polynilpotent
Multipliers of Finitely Generated Abelian Groups, {\it
International Jornal of Mathematics, Game Theory and Algebra},
{\bf 16:1}, (2006) 93-102.
\bibitem{8} M. R. R. Moghaddam, The Baer invariant of a Direct
Product, {\it Archiv der Math.}, {\bf 33}, (1979) 504-511.
\bibitem{9} M. R. R. Mogaddam and B. Mashayekhy, Higher Schur
Multiplicator of a Finite Abelian Group, {\it Algebra Colloquium
}, {\bf 4:3}, (1997) 317-322.
\bibitem{10} M. R. R. Moghaddam, B. Mashayekhy, and S. Kayvanfar,
 The Higher Schur Multiplicator of Certain Class of Groups,
{\it Southeast Asian Bulletin of Mathematics}, {\bf 27}, (2003)
121-128.
\bibitem{11} H. Neuman, {\it Varieties of Groups},
Berlin-Heidelberg-New York, Springer, 1967.
\bibitem{12} M. Parvizi and B. Mashayekhy , On Polynilpotent
Multipliers of Free Nilpotent Groups, {\it Communication in
Algebra}, {\bf 34:6}, (2006) 2287-2294.
\bibitem{13} R. R. Struik, On Nilpotent Products of Cyclic Groups,
{\it Canada. J. Math.}, {\bf 12}, (1960) 447-462.
\end{thebibliography}
\end{document}